\documentclass{article}
\usepackage{amssymb,amsmath,eucal}

\input{epsf}

\def\bfh{\mathbf h}

\def\frA{\mathfrak A}

\def\frS{\mathfrak S}

\def\frp{\mathfrak p}
\def\frq{\mathfrak q}
\def\frr{\mathfrak r}

\def\fru{\mathfrak u}

\def\cK{\mathcal K}

\def\Q {{\mathbb Q }}
\def\R {{\mathbb R }}
\def\C {{\Bbb C }}
\def\K{{\Bbb K}}

\def\Z{{\Bbb Z}}

\def\la{\langle}
\def\ra{\rangle}

\newcommand{\dom}{\mathop{\rm dom}\nolimits}
\newcommand{\im}{\mathop{\rm im}\nolimits}

\def\GL{{\rm GL}}
\def\SL{{\rm SL}}

\def\U{{\rm U}}

\def\cP{\mathcal P}
\def\cM{\mathcal M}

\def\ov{\overline}
\def\phi{\varphi}
\def\epsilon{\varepsilon}
\def\kappa{\varkappa}
\def\le{\leqslant}
\def\ge{\geqslant}

\def\wh{\widehat}

\newcounter{sec}

\newcounter{punct}[sec]

\def\punct{\refstepcounter{punct}{\arabic{sec}.\arabic{punct}.  }}

\newtheorem{theorem}{Theorem}[sec]
\newtheorem{proposition}[theorem]{Proposition}

\newtheorem{lemma}[theorem]{Lemma}

\def\G{\mathbb{G}}

\def\N{\mathbb{N}}

\def\COUNTERS{\addtocounter{sec}{1}
              \setcounter{punct}{0}
          \setcounter{equation}{0}
          \setcounter{theorem}{0}
          \setcounter{problem}{0}
          }
          
          \def\sm{\smallskip}
          \def\ov{\overline}
          \def \wt{\widetilde}
         \def \cO{\mathcal{O}}
\def\1{\mathbf 1}

\def\Gal{\mathrm {Gal}}

\def\whi{\mathrm{white}}
\def\bla{\mathrm{black}}
\def\blu{\mathrm{blue}}
\def\red{\mathrm{red}}
\def\yel{\mathrm{yellow}}

\def\white{\mathrm{\whi}}
\def\black{\mathrm{\bla}}
\def\blue{\mathrm{\blu}}
\def\yellow{\mathrm{\yel}}

\def\cP{\mathcal P}
\def\cQ{\mathcal Q}
\def\cR{\mathcal R}

\def\diag{\mathrm{diag}}

\def\G{\mathbb{G}}
\def\K{\mathbb{K}}
          
\begin{document}

\begin{center}
{\bf\Large
Symmetric groups and

\medskip

checker triangulated surfaces}

\bigskip


\sc \large 
Yury A. Neretin%
\footnote{Supported by the grant FWF, P28421.}
\end{center}

{\small We consider triangulations of surfaces with edges painted 
three colors so that edges of each triangle have different
colors. Such structures arise as Belyi data (or Grothendieck dessins d'enfant),
on the other hand they enumerate pairs of permutations determined up to a common conjugation.
The topic of these notes is links of such combinatorial
structures with infinite symmetric groups
and their representations.
}

\vspace{22pt}

{\sc 

1. Checker triangulations and Belyi data.

2. Correspondence with symmetric groups.

3. Category of checker cobordisms.

4. The infinite tri-symmetric group.

5. Spherical functions.

6. Convolution algebras and their limits.

7. A pre-limit convolution algebra.}

\section{Checker triangulations and Belyi data}

\COUNTERS

\begin{figure}
 $$\epsfbox{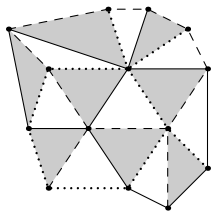}$$
 \caption{A piece of a checker surface. We draw full lines, dashed-lines, dotted-lines
 instead of colors%
 \label{fig:piece}}
\end{figure}

{\bf \punct Checker triangulations.} Consider an oriented compact closed two-dimensional surface 
$P$ (we do not require that a surface  be connected). Consider a finite graph $\gamma$  on $P$ separating $P$ into
triangles.
 Let  edges of the graph be colored  {\it blue, red, or yellow} and let
each triangle  have edges of all 3 types. We say that a triangle is {\it white} if 
the order of
blue, red, yellow edges is clockwise  and {\it black} in the opposite case, 
see Fig.~\ref{fig:piece}. We call such structures $(P,\gamma)$  {\it checker surfaces}.

Two checker surfaces $(P,\gamma)$ and $(P,\gamma')$ are {\it equivalent} if 
there is an orientation preserving homeomorphism $P\to P'$ identifying colored graphs $\gamma$ and $\gamma'$.
Therefore checker surfaces are purely combinatorial objects, we also can regard them as
simplicial cell complexes with colorings of triangles and edges.

We denote the set of all checker  surfaces having $2n$ triangles by $\Xi_n$, 

\sm

\begin{figure}
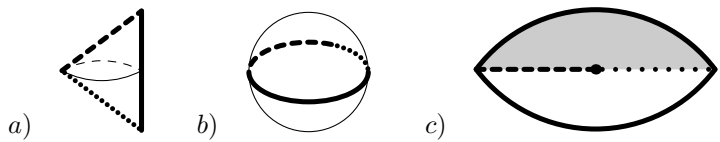

  $$a) \quad \epsfbox{checker.3}\qquad b) \quad \epsfbox{checker.4}\qquad c) \quad\epsfbox{checker.5 } $$
  \caption{a) A double triangle.
  \newline
b)  A double triangle drawn as a sphere.
\newline
c) This position of two triangles is admissible.%
\label{fig:double}
}
\end{figure}

{\sc Example.} Consider a sphere with a  cycle composed of blue, red, and yellow edges,
see Fig.~\ref{fig:double}a. 
We call such checker surface  a {\it double triangle}.

\sm

{\sc Remarks.} a) Under our conditions orders of all vertices are even.
Edges adjacent to a given vertex have two colors, and these colors are interlacing,
see Fig.~\ref{fig:piece}.
Therefore we can split vertices 
 into 3 classes, blue, red, yellow:
we say that a {\it vertex is blue}, if opposite edges are blue, etc.
(colors of vertices are
not shown on  Fig.~\ref{fig:piece}).

\sm

b) An intersection of two white  (or two black) triangles is the empty set or a vertex.
An intersection
of a black and a white triangles can be the empty set, a vertex, an edge, two edges, three edges, see 
Fig.~\ref{fig:double}c,
\ref{fig:double}a.
We avoid a term triangulation since the most common definition of  triangulations forbids the last two possibilities.
Formally, our objects satisfy a definition of a simplicial cell complex, see, e.g., \cite{Hat}. 
\hfill $\boxtimes$

\sm

c) We admit empty surfaces.

\sm

{\sc Remark.}
Apparently, such colored structures  firstly appeared in combinatorial topology
in \cite{Pez}, see also \cite{FGG}, \cite{Gai}.
\hfill $\boxtimes$

\sm

{\bf\punct Dessins d'enfant.}
There is the following equivalent language for description of the same objects. 
We can remove red and yellow edges and blue vertices. After this we get a  graph with blue edges
on the surface, vertices are painted red and yellow, and colors of the ends of  each edge 
are different. Therefore we get a bipartite graph on the surface, a coloring of vertices is determined by the graph
modulo a reversion of colors of all vertices.
Such a graph separates the surface into polygonal domains.

\begin{figure}
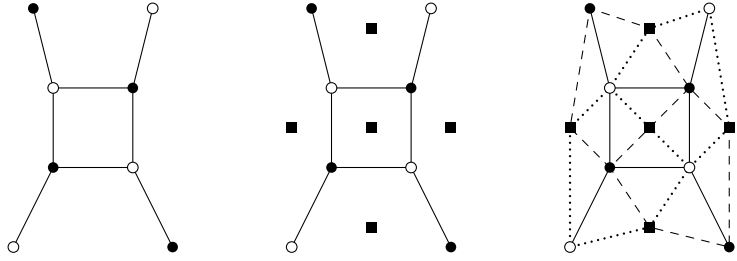

$$
\epsfbox{checker.6} \qquad  \qquad \epsfbox{checker.7} \qquad \qquad  \epsfbox{checker.8}
$$
 \caption{A reconstruction of checker structure from a dessin d'enfant.%
 \label{fig:dessin}}
\end{figure}

 To reconstruct 
a checker structure, we put a point  to interior  of each domain, connect it by non-intersecting arcs with
vertices of the polygon, and paint these arcs in a unique admissible way,
see Fig.~\ref{fig:dessin}.

\sm

{\bf\punct Belyi data.}
Denote by $\wh\C$ the Riemann sphere, by $\ov \Q$  the algebraic closure of $\Q$,
 by $\Gal(\ov\Q/\Q)$ the Galois group of $\ov \Q$ over $\Q$.
 Recall  the famous
{\it Belyi theorem} \cite{Bel-1}, \cite{Bel-2}, see expositions in \cite{LZ}, \cite{G-G}.

\begin{theorem}
	Let $P$ be a nonsingular complex algebraic (projective) curve.
A covering holomorphic map $\pi:P\to\ov\C$ whose  critical values 
are contained in the set $\{0,1,\infty\}$
 exists if and only if $P$ is defined over  $\ov \Q$. 
 For a curve defined over $\ov\Q$ a map $\pi$ can be chosen defined over $\ov\Q$.
\end{theorem}

Words a {\it a curve $P$ is defined over $\ov\Q$} mean that $P$ can be determined
by a system of algebraic equations
in a complex projective space $\C\mathbb{P}^N$ with coefficients in $\ov \Q$.

Such functions  are called {\it Belyi functions}, a pair $(P,\pi)$ is called {\it Belyi data}.

The Galois group $\Gal(\ov \Q/\Q)$ acts on Belyi data changing coefficients
of equations determining  curves.

\sm

{\bf \punct Correspondence between Belyi data and checker surfaces.}
The real line $\R$ splits $\wh \C$ into two half-planes. Let us paint
the segment $[\infty,0]$ blue, the segment $[0,1]$ red, and the segment $[1,\infty]$ yellow.
The upper half-plane of $\wh \C$ becomes white and the low half-plane becomes black. 
Thus we get a structure of a double triangle on $\wh \C$, see Fig.~\ref{fig:double}b. 

\sm

For a Belyi map $P\to \wh\C$ we take the preimage of $\R$. It is a colored
graph on $\R$ splitting $P$ into white and black triangles. Thus we get
a structure of a checker surface on $P$.

\sm

Conversely, let us for a checker surface $(P,\gamma)$
define a complex structure on $P$ and a holomorphic function $P\to\wh \C$.
 Consider 
an equilateral triangle $T^{\white}$ whose sides are colored blue, red, yellow 
clockwise and an equilateral triangle $T^{\black}$ whose sides are colored blue, red, yellow 
anti-clockwise. We identify $T^{\white}$ conformally with the white half-plane of $\wh\C$ sending vertices to
0, 1, $\infty$ and edges to segments $[-\infty,0]$, $[0,1]$, $[1,\infty]$ according colors
(such a map exists and is unique).
 We also identify $T^{\black}$ with the black half-plane in a similar way.
We can think that a surface $(P,\gamma)$ is glued from copies of $T^{\white}$ and $T^{\black}$. This defines 
  a complex structure on $P$ outside vertices. If we have $2k$ triangles adjacent to a vertex $v$,
 then
 a  punctured  neighborhood of $v$ has a structure of $k$-covering of a punctured circle $\dot D:
 \,0<|z|<\epsilon$. This covering is itself holomorphically equivalent to $\dot D$,
 and this gives us a chart on $P$ in a neighborhood of $v$. Thus we get a complex structure on
 $P$.

Next, we take a map from $P$ to $\wh\C$ sending   white triangles to the white half-plane and 
black triangles to the black  half-plane. By the Riemann--Schwarz reflection principle the map is
holomorphic at interior points of edges. By the theorem on removable singularities
the map is also holomorphic at vertices.

Thus we get Belyi data. It is known, that combinatorial data determine a Belyi map
up to a natural equivalence (upto an automorphism of $P$). 
Notice, that in our case $P$ can be disconnected, so we have Belyi functions
on all components of the curve.

\bigskip

In particular the  group $\Gal(\ov \Q/\Q)$ acts on the space of checker surfaces. In 1984 Grothendieck
\cite{Gro}
 initiated a program of investigation of $\Gal(\ov \Q/\Q)$ based on dessins d'enfant, see e.g.,  \cite{Sch}, \cite{S-Loshak},
\cite{ShV}, \cite{LZ}, \cite{G-G}, \cite{SV}, \cite{Shabat}.

\section{Correspondence with symmetric groups}

\COUNTERS

{\bf\punct The basic construction.}
Consider the product $G_n=S_n\times S_n\times S_n$ of 3 copies of the  symmetric group $S_n$.
An element of this group is a triple of permutations, we denote such triples as 
\begin{equation}
p=\bigl(p^{\blue}, p^\red, p^\yellow\bigr).
\label{eq:wt-g}
\end{equation}
Denote by $K_n$ the diagonal subgroup in $G_n$, it consists of triples of the form
\begin{equation}
\bfh=(h,h,h).
\label{eq:wt-h}
\end{equation}

Take a triangle $T^{\white}$, paint its sides  blue, red, yellow clockwise. Consider $n$ copies 
of the triangle $T^{\white}$, paint them  white and attribute to the triangles labels 1, 2, \dots, $n$.
Thus we get triangles $T^{\white}_1$, \dots, $T_n^{\white}$.

\begin{figure}
 $$\epsfbox{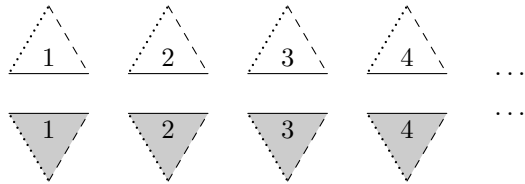}$$
 \caption{The collection of triangles $T_j^{\white}$, $T_j^{\black}$.%
 \label{fig:triangle}}
\end{figure}

Consider another triangle $T^{\black}$ obtained from $T^{\white}$ by a reflection. Its sides are colored
 blue, red, yellow anti-clockwise. Take $n$ copies of this triangle $T^{\black}_1$, \dots, $T_n^{\black}$,
 see Fig.~\ref{fig:triangle}.
 
 Fix an element $p\in G_n$,
For each $j=1$,\dots, $n$ we glue the triangle $T_j^\whi$ with
the triangle  $T_{p^\blu j}^\bla$ along blue sides
according orientations. Repeat the same operation for $p^\red$ and $p^\yel$.
We get a two-dimensional surface equipped
with a checker structure. We also have an additional structure, namely we have a
bijective map ({\it labeling}) from the set $\{1,\dots,n\}$ to the set of white triangles 
and the bijective map from the same set to the set of black triangles. 

We call the object $\cP$ obtained in this way a {\it completely labeled checker surface}. Denote by
$\wt \Xi_n$ the set of all  such surfaces.

\sm

{\bf\punct The inverse construction.} Let us describe the inverse map $\wt\Xi_n\to G_n$.
Consider a completely labeled surface $\cP$. Take a blue edge $v$. Let $i_v$ be the label on the white side of $v$,
and $j_v$ the label on the black side. Then we set
$$
p^{\blu}:i_v\mapsto j_v.
$$
and get $p^\blu$. In the same way we reconstruct $p^\red$ and $p^\yellow$.

\sm

{\bf \punct Remarks.}
a) Thus for any element of $G_n$ we get a surface with $2n$ triangles and $3n$ edges. 

\sm

b) Blue vertices of $\cP$ are in one-to-one
correspondence with independent cycles of the permutation $(p^\yel)^{-1} p^\red $. 
The order of a vertex is the duplicated length of a cycle.
A similar statement holds for red and yellow vertices.

\sm

c) Take  the subgroup $U_p\subset S_n$ generated by $(p^\yel)^{-1} p^\blue$ and  $(p^\yel)^{-1} p^\red $.
Components $\cP_j$ of $\cP$ are in one-to-one correspondence with orbits $\cO_j$ of $U_p$ on $\{1,\dots,n\}$.
A given orbit splits into independent cycles of  
 $(p^\yel)^{-1} p^\red $, this gives us an enumeration of blue vertices of the component.
 So we have an expression for the Euler characteristic of a component in group-theoretical terms:
 \begin{multline*}
 \chi(\cP_j)=-\# (\cO_j)+
 \begin{Bmatrix}
 \text{number of cycles of}\\
 \text{$(p^\yel)^{-1} p^\blue$ on $\cO_j$}
 \end{Bmatrix}
+\\+
 \begin{Bmatrix}
\text{number of cycles of}\\
\text{of $(p^\yel)^{-1} p^\red$ on $\cO_j$}
\end{Bmatrix}
+
 \begin{Bmatrix}
\text{number of cycles of}\\
\text{$(g^\red)^{-1} g^\blue$ on $\cO_j$}
\end{Bmatrix}.
  \end{multline*}
 
 
 \sm

 d) The pass $ p\mapsto p^{-1}$ is a reversion of colors of triangles black$\leftrightarrows$white
 and a reversion of orientations  of all components.
 
 \sm
 
 e) For any element $\bfh\in K_n$ the corresponding surface is a disjoint union of double triangles.
 
 \sm

 \sm
 
 {\bf \punct Concatenations.%
 \label{ss:concatenation}} Now let us describe the product in the group $G_n$
 on the geometric language.
 
 Let
 $p$, $ q \in G_n$. Consider the corresponding elements $\cP$, $\cQ\in \wt\Xi_n$.
 For each $j=1$, 2, \dots, $n$ we identify $j$-th black triangles of $\cQ$ and $j$-th white triangle of
 $\cP$ according colors of edges. In this way we get a 2-dimensional simplicial cell complex $\cR^{\circ\circ}$ having
 $3n$ two-dimensional cells, namely:
 
 \sm
 
 --- $n$ labeled images of white triangles of $\cQ$;
 
 \sm
 
 --- $n$ labeled images
 of black triangles of $\cP$;
 
 \sm
 
 --- $n$ (labeled) triangles that are results of gluing (let us call such cells
 {\it grey}).
 
 \sm
 
 Each edge  has a color (blue, red, or yellow) and is contained in 3 triangles (white, black, and grey).
 Next, we remove interiors of all grey triangles from the simplicial cell complex $\cR^{\circ\circ}$ (and forget their labels).
 We get a complex $\cR^\circ$ having $n$ white triangles and $n$ black triangles.
 These cells inherit labelings. Moreover, each
 edge is contained in precisely two triangles. 
 
\begin{figure}
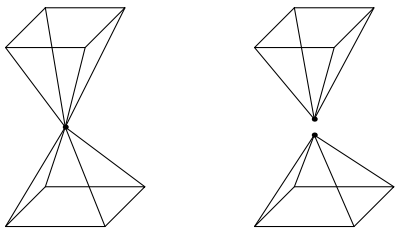

 $$\epsfbox{checker.9} \qquad \qquad \epsfbox{checker.10}$$
 \caption{A normalization.%
\label{fig:normalization}}
\end{figure}

 In fact, we get a surface, but some vertices of the surface are glued one with another, see
 Fig.~\ref{fig:normalization}.a.
 Cutting such  gluings%
 \footnote{
 Let us describe this cutting (a {\it normalization})
 formally (this is a special case of procedure of normalization of pseudomanifolds
 from \cite{GorMPh}). Let us equip the  space
 $\cR^\circ$ by some natural metric (for instance we assume that all triangles are equilateral 
triangles with side 1, this determines a 'geodesic' distance on the whole complex. 
 Let us remove all vertices, consider the geodesic distance on this space, and complete our space
 with respect to the new distance.}
 we get a surface $\cR$ corresponding  to the product
 $r= p q\in G_n$.
 
\sm

{\sc Remark.}
Of course we can transpose an order of the  gluing and the making  holes. Namely, let $\cP$ and $\cQ$ be
completely
labeled checker surfaces. We remove interiors of black triangles from $\cP$ and interiors of white triangles
from $\cQ$ and get new cell complexes $\cP^\square$ and $\cQ^\square$.
For each $j\le n$ we identify the boundary of former $j$-th black triangle of $\cQ$ and the boundary
of former $j$-th white triangle of $\cP$ according colors of edges. We get a two-dimensional cell complex and normalize it
as above.
\hfill $\boxtimes$

 \sm

{\bf \punct Notation. Conjugacy classes and double cosets.}
 Let $H$ be a group, $L$, $M\subset H$  subgroups. We denote
 by $L\setminus H/M$ the space of {\it double cosets} of $H$ with respect to 
 $L$, $M$, i.e., the quotient of $H$ with respect to the equivalence
 $$
 h\sim l h m,\qquad\text{where $l\in L$, $m\in M$.}
 $$
 An element of this space is a subset in $H$ of the form $L\cdot h\cdot M$.
 
 \sm
 
 Denote by $H/\!\!/L$ the space of conjugacy classes  of a group
  $H$ with respect to a subgroup $L$, i.e.,
 the quotient of $H$ with respect to the equivalence
 $$
 h\sim l^{-1} h l, \qquad\text{where $l\in L$.}
 $$

 {\sc Remark.}
Conjugacy classes can be regarded as a special case of double cosets for the following reason.
Consider the group $H\times L$ and the embedding $\iota:L\to H\times L$ given by $\iota(l)= (l,l)$.
Obviously,
$$
H/\!\!/L\,\,\simeq\,\, \iota(L)\setminus (H\times L)/\iota(L).
$$

\sm

{\bf \punct Relabeling.} Let us multiply an element $p\in G_n$ by an element $\bfh\in K_n$, see
 (\ref{eq:wt-g}) and (\ref{eq:wt-h}). Clearly the operation 
 $p\mapsto \bfh  p$
 is equivalent to a permutation of 
 labels on black triangles of the corresponding surface.
 A right multiplication $p\mapsto  p\bfh^{-1}$
 is equivalent to a permutation of labels on white triangles.
 
 \sm
 
 Consider the space of  double cosets of $K_n\setminus G_n/K_n$ of $G_n$ with respect to $K_n$.
 This means that we consider elements of $\wt \Xi_n$ up to a permutation of labels.
 In other words, we forget labels.
 
 Thus we get a bijection 
 $$K_n\setminus G_n/K_n\, \longleftrightarrow \,\Xi_n.$$

 \sm
 
 This also gives a geometrical  description of pairs of permutations determined up to a common conjugation, i.e., upto the equivalence
 $$
 (g_1,g_2)\sim (h^{-1} g_1 h, h^{-1} g_2 h).
 $$
 Indeed, we have the space of conjugacy classes of $S_n\times S_n$ with respect to the diagonal
 subgroup $\diag(S_n)$.
 According the remark in previous subsection, 
 $$(S_n\times S_n)/\!\!/\diag(S_n)\simeq
\diag(S_n)\setminus (S_n\times S_n\times S_n)/\diag(S_n)=
  K_n\setminus G_n/K_n.$$

\section{Category of checker cobordisms}

\COUNTERS

{\bf \punct Infinite symmetric group.}
Denote by $S_\infty$ the group of finitely supported permutations of $\N$,
by $\ov S_\infty$ the group of all permutations. The group
$S_\infty$ is a countable discrete group, the group $\ov S_\infty$
has  cardinality continuum and 
is equipped with a natural topology discussed  in the next section.
As usual, we can represent elements of these groups by 
0-1-matrices having precisely one 1 in each column and each row.

Consider the group
$$G:=S_\infty\times S_\infty\times S_\infty$$
and its diagonal $K\simeq S_\infty$.
We  repeat the construction of the previous section for $n=\infty$.
For any $g\in G$ we get a countable disjoint union of 
compact checker surfaces, all but a finite number of components are double triangles,  moreover
for all but a finite number of components labels on the black and white sides 
of a double triangles coincide.
As above, we get a correspondence between the group and the set of completely labeled checker surfaces.

\sm 

{\bf\punct Multiplication of double cosets.%
\label{ss:multiplication-cosets}}
Denote by $\Z_+$ the set of non-negative integers.
For $\alpha\in \Z_+$ denote by $K[\alpha]\subset K$ the subgroup in $K=S_\infty$
consisting of permutations fixing points $1$, \dots, $\alpha\in \N$; for $\alpha=0$ we set
$K[0]:=K$. The subgroups $K[\alpha]$ are isomorphic to $S_\infty$.

Consider double cosets spaces
$$\cM(\alpha,\beta):=K[\alpha]\setminus G/K[\beta].$$
It turns out to be that for any $\alpha$, $\beta$, $\gamma\in \Z_+$ there is a natural  multiplication
$$
\cM(\alpha,\beta) \times \cM(\beta,\gamma)\,\to\,\cM (\alpha,\gamma)
$$
defined in the following way. 
Consider a sequence $\theta_j[\beta]$ in $S_\infty$
given  by 
$$
\theta_j[\beta]:=
\begin{pmatrix}
 1_\beta&0&0&0\\
 0&0&1_j&0\\
 0&1_j&0&0\\
 0&0&0&1_\infty
\end{pmatrix}.
$$
By $\Theta_j[\beta]$ denote the corresponding element of $K$:
$$ \Theta_j[\beta]=\bigl(\theta_j[\beta],\theta_j[\beta],\theta_j[\beta]\bigr)   \in K[\beta].$$
Let $\frp\in \cM(\alpha,\beta)$, $\frq\in \cM(\beta,\gamma)$. Choose   their representatives $p\in \frp$, $q\in \frq$.
The following statements are semi-obvious if to look to them for a sufficiently long time
(see formal proofs in \cite{GN}).

\begin{lemma}
  The following sequence of double cosets
 \begin{equation}
 K[\alpha] \cdot p\cdot \Theta_j[\beta]\cdot q\cdot K[\gamma]\in \cM(\alpha,\gamma)
 \label{eq:const}
 \end{equation}
 is eventually constant. Its limit%
 \footnote{I.e., a stable value. However,  $\cM(\alpha,\gamma)$
 is a quotient of a discrete topological space, so it has a discrete topology, and we
 have a limit in the sense of the formal definition of a limit.}
 does not depend on a choice of representatives
 $p\in \frp$, $q\in \frq$.
\end{lemma}

Denote by $\frp\circledast \frq$ the limit of the sequence (\ref{eq:const}).
The multiplication $\circledast$ is associative in the following sense:

\begin{lemma}
 For any $\alpha$, $\beta$, $\gamma$, $\delta\in \Z_+$ and any 
 $$\frp\in \cM(\alpha,\beta), \quad\frq\in\cM(\beta,\gamma),
 \quad\frr\in\cM(\gamma,\delta),
 $$
 we have 
 $$
( \frp \circledast \frq) \circledast \frr=  \frp \circledast (\frq \circledast \frr).
 $$
\end{lemma}

Thus we can define the following category $\cK$. The set of its objects is $\Z_+$ and
the set of morphisms from $\beta$ to $\alpha$ is $\cM(\alpha,\beta)$.

This operation is a representative of huge zoo of train constructions for infinite dimensional groups,
see end of the present section.
The main property of $\circledast$-multiplication  is Theorem \ref{th:multiplicativity} below.

\sm

{\bf\punct The involution on $\cK$.} The map $p\mapsto p^{-1}$ induces  maps
$$
\cM(\alpha,\beta)\to\cM(\beta,\alpha),
$$
we denote these maps by 
$\frp\mapsto \frp^*$. Obviously, $*$ satisfies the usual properties of involutions
$$
(\frp\circledast \frq)^*=\frq^*\circledast \frp^*, \qquad \frp^{**}=\frp.
$$

\sm

{\bf\punct Explicit description of the product of double cosets.\label{ss:description-cosets}}
We say that a
{\it labeled checker surface} is the following collection of combinatorial data:

\sm

--- a finite disjoint union $P$ of checker surfaces;

\sm

--- an injective map from the set $\{1, \dots, \alpha\}$ to the set of black triangles
of $P$ and 
an  injective map from the set $\{1, \dots, \beta\}$ to the set of white triangles
(we call these maps by {\it labelings}).

\sm

--- additionally, we require that any double-triangle component of $P$ 
keeps at least one label.

\sm

Denote by $\wt\Xi[\alpha,\beta]$ the set of such objects.

Clearly, there is a canonical one-to-one correspondence 
$$\cM(\alpha,\beta)\, \longleftrightarrow\, \wt\Xi[\alpha,\beta].$$
 Namely, we take a representative of a double coset $\frp$ and
draw the corresponding completely labeled checker surface. A pass to the double coset means forgetting
white labels with numbers $>\beta$ and black labels $>\alpha$. Finally, we remove all double triangles
without labels and get a finite object%
\footnote{After the previous operation we have a countable collection of label-less double triangles,
which do not keep any information. It is more convenient to remove them.}.

\sm

The product of double cosets corresponds  to the concatenation
(as in Subsect. \ref{ss:concatenation}). Namely, let $\cP$, $\cQ$ be labeled checker surfaces
corresponding to $\frp$ and $\frq$. For each $j=1$, \dots, $\beta$
we remove interiors of labeled black triangles
of $\cQ$ and  labeled  white triangles of $\cP$
(so we get $\beta$ 'holes' on each surface). For each $j\le\beta$ we identify boundary
of the former black triangle of $\cQ$ and the former white triangle of $\cP$ according colors of edges.
We get a new 2-dimensional simplicial cell complex and normalize it as above. This surface
inherits $\alpha$ labels on white triangles and $\gamma$ labels on black triangles.
Finally, we remove label-less double triangles.

\sm

Thus  we get an operation similar to a concatenation of cobordisms, we glue 
cell complexes with additional structure instead of manifolds with boundaries,
also our operation includes a normalization.

\sm

{\sc Remark.}
The operation of normalization does not arise if labeled black triangles of 
$\cQ$ (or labeled white triangles of $\cP$) have no common vertices.
\hfill $\boxtimes$

\sm

{\sc Remark.} The $\circledast$-product on $\wt\Xi[0,0]$ is a disjoint union of
label-less
checker surfaces,
$$
\cP \circledast\cQ=\cP\coprod\cQ.
$$
In particular, this operation is commutative.
\hfill $\boxtimes$

\sm

{\bf \punct Multiplicativity theorem.}
Let $\rho$ be a unitary representation of the group $G=S_\infty\times S_\infty\times S_\infty$
in a Hilbert space%
\footnote{For many classes of groups (as finite groups, compact groups, semisimple Lie groups,
nilpotent Lie groups, semisimple $p$-adic groups)
any irreducible representation of $G_1\times G_2$ is a tensor product
of an irreducible representation $\rho_1$ of $G_1$ and an irreducible
representation $\rho_2$ of $G_2$.
In fact, this holds if at least one of groups has type I, see, e.g., [18], \S5,
\S 13.1. For the group $S_\infty$ this implication does not valid.
Also, in the next section we observe that actually we work with a certain completion of 
the group $G$ and this completion is not a product of groups.}
$V$. For each $\alpha\in \Z_+$ denote by $V[\alpha]$ the space of 
all $K[\alpha]$-fixed vectors  $v\in V$,
$$
\rho(\bfh) v=v \qquad\text{for all $\bfh\in K[\alpha]$}.
$$

{\sc Remark.}
Obviously,
$$
V[0]\subset V[1]\subset V[2]\subset \dots \subset V.
$$
In this moment we do not require  $\cup_\alpha V[\alpha]$ be dense in $V$.
However, below we will observe that this requirement is natural.
\hfill $\boxtimes$

\sm

Denote by $P[\alpha]$ the orthogonal projection operator $V\to V[\alpha]$.

For $\alpha$, $\beta\in \Z_+$, and  $p\in G$ consider the operator
$$
\ov \rho (p): V[\beta]\to V[\alpha]
$$
given by
$$
\ov\rho_{\alpha,\beta}(p):=P[\alpha] \rho(p)\Bigr|_{V[\beta]}.
$$
In other words, we write $\rho(p)$ as a block operator
$$
\rho_{\alpha,\beta}(p)=\begin{pmatrix}
         R_{11}&R_{12}\\R_{21}&R_{22}
        \end{pmatrix}:\quad V[\beta]\oplus V[\beta]^\bot \to V[\alpha]\oplus V[\alpha]^\bot
$$
and set 
$$
\ov \rho(p):=R_{11}.
$$
Then $\ov \rho(p)$ depends only on a double coset containing $p$
(this is a standard fact), so actually we get a function
$\ov\rho_{\alpha,\beta}(\frp)$ on $\cM(\alpha,\beta)$.

\begin{theorem}
	\label{th:multiplicativity}
 The operators $\ov\rho_{\cdot,\cdot}(\cdot)$ determine a representation of the category $\cK$ of double cosets,
 i.e., for any $\alpha$, $\beta$, $\gamma\in \Z_+$ for any $\frp\in \cM(\alpha,\beta)$, $\frq\in \cM(\beta,\gamma)$,
 $$
 \ov\rho_{\alpha,\beta}(\frp) \,\ov\rho_{\beta,\gamma}(\frq)=\ov\rho_{\alpha,\gamma}\bigl(\frp\circledast\frq\bigr)
 .
 $$
 Moreover, this is a $*$-representation of the category $\cK$, i.e.,
 $$
 \ov\rho_{\alpha,\beta}(\frp)^*=\ov\rho_{\beta,\alpha}(\frp^*). 
 $$
\end{theorem}

See \cite{Ner-imrn}, \cite{Ner-umn}.

Thus  any unitary representation $\rho$ of $G$ generates
a $*$-representation $\ov\rho$ of the category $\cK$. However, generally speaking
$\rho$ can not be reconstructed from
$\ov\rho$. In the next section we describe a more perfect version of this correspondence.

\sm

{\bf \punct Zoo.}
Apparently, a first example of a multiplication on double cosets spaces
$K\setminus G/K$ and a multiplicativity theorem was discovered
by R.~S.~Ismagilov \cite{Ism} for a pair $G\supset K$, where $G=\SL_2$
over a non-locally compact non-Archimedean field 
and $K$ is $\SL_2$ over integer elements of the field.

Later many cases of multiplication were examined by G.~I.~Olshanski,
see \cite{Olsh-semigr}. In particular, he considered the case $G=S_\infty\times S_\infty$
and $K=K[\alpha]$, see \cite{Olsh-chip}; 
there is a well-developed theory of representations of this group,
see \cite{KOV}, \cite{BO} and further references in this book.

In \cite{Ner-book}, Sect. VIII.5, it was observed that multiplicativity theorems  are a quite
general phenomenon,  the main obstacle for a further progress was a problem of explicit descriptions of double coset spaces.
The construction of this section was obtained in \cite{Ner-imrn}.
It has numerous variations with different pairs of groups $G\supset K$
related to infinite symmetric groups, see \cite{Ner-umn}, \cite{GN} as
\begin{align*}
&G= (S_\infty)^n, \qquad K=\mathrm{diag} (S_\infty);
\\
&G= S_{n\infty},\qquad K= S_\infty \ltimes (\Z_n)^\infty;
\\
&G= S_{n\infty},\qquad K= S_\infty \ltimes (S_n)^\infty;
\\
&G=S_{n\infty}, \qquad K= (S_\infty)^n.
\end{align*}
For the last case Nessonov (see \cite{Nes}, \cite{Ner-umn})
obtained a classification of all $K$-spherical
representations of $G$.

On the other hand, there is parallel  picture for infinite dimensional real classical
groups, \cite{Olsh-GB}, \cite{Ner-book} (Sect. IX.3), \cite{Ner-class-1}, \cite{Ner-class-2}.
There is a poorly  studied mass construction of the same type for infinite-dimensional groups
over $p$-adic field \cite{Ner-p-adic} and an example over finite fields
\cite{Ner-finite}. On the other hand, this story has a counterpart
for groups of transformations of measure spaces \cite{Ner-book}, Chapter X, \cite{Ner-jfa}.

\section{The infinite tri-symmetric group}

\COUNTERS

Here we describe a natural completion $\G$ of $G=S_\infty\times S_\infty\times S_\infty$,
$$
S_\infty\times S_\infty\times S_\infty \,\, \subset\,\, \G\,\,\subset\,\, \ov S_\infty\times \ov S_\infty\times \ov S_\infty
$$
and discuss representations of this completion.

\sm

{\bf\punct Representations of a full symmetric group $\ov S_\infty$.}
We say that a sequence $g_j\in \ov S_\infty$ converges to $g$ if for each $k\in\N$
we have $g_j k=g k$ for sufficiently large $j$. This determines a structure of a
complete separable 
 topological group on $\ov S_\infty$ (and this topology is unique in any reasonable sense, see \cite{KR}).
 
 Denote by $\ov S_\infty[\alpha]\subset \ov S_\infty$ the subgroup consisting of elements fixing
 points 1, 2, \dots,  $\alpha\in\N$. These subgroups form a base of neighborhoods of unit in $S_\infty$
 (and this  can be regarded as another definition of topology in $S_\infty$).

 A classification of unitary representations of $\ov S_\infty$ was obtained by A.~Lieberman \cite{Lie}
 (see also expositions in \cite{Olsh-kiado}, \cite{Ner-book}, Sect. VIII.1-2),
 the classification  is relatively simple but it is not necessary for us in  the sequel%
 \footnote{Any irreducible unitary representation of $\ov S_\infty$ has the following
 form. Fix $\alpha\in\Z_+$ and an irreducible representation
 $\rho$ of $S_\alpha$.  Consider a subgroup $S_\alpha\times \ov S_\infty[\alpha]$ and its representation
 that is trivial on $\ov S_\infty[\alpha]$ and equal $\rho$ on $S_\alpha$. 
 The induced representation of $\ov S_\infty$ is irreducible, and all irreducible representations
 of $\ov S_\infty$ have this form. Any unitary representation is a direct sum of irreducible representations.
 A decomposition of a unitary representation into a direct sum of irreducible
 representations is unique.}.
 Representations of $S_\infty$ that can be extended
 to $\ov S_\infty$ can be described in the following form without explicit reference to the topology.

  \begin{proposition}
  Let $\rho$ be a unitary representation of $S_\infty$ in a Hilbert space $V$. As above, denote by
  $V[\alpha]$ the space of $S_\infty[\alpha]$-fixed vectors. 
  
  {\rm a)}
  A representation $\rho$ admits
  a continuous extension to $\ov S_\infty$ if and only if $\cup_\alpha V[\alpha]$ is dense in $V$.
  
  \sm
  
  {\rm b)} An irreducible representation $\rho$ admits a continuous extension to $\ov S_\infty$ if 
  and only if $\cup_\alpha V[\alpha]$ is non-zero.
 \end{proposition}

 
 {\bf \punct Tri-symmetric group.} Consider the group
 $\ov S_\infty\times\ov S_\infty\times \ov S_\infty$ and its diagonal
 subgroup $\K\simeq\ov S_\infty$. We define the {\it tri-symmetric group} $\G$ as
 the subgroup in  $\ov S_\infty\times\ov S_\infty\times \ov S_\infty$
 generated by $\K$ and $G=S_\infty\times S_\infty\times S_\infty$.
 In other  words
 $$
 (p^\blu, p^\red, p^\yel)\in \G\qquad \Longleftrightarrow\qquad
 p^\blu(p^\red)^{-1},\,\,  p^\red(p^\yel)^{-1}\in S_\infty
 .
 $$
 The coset space $\G/\K$ is countable. We define a topology on $\G$ assuming that $\K$ has the natural topology
 and the space $\G/\K$ is discrete.
 
 For the group $\G$ we repeat the considerations
 of Subsections \ref{ss:multiplication-cosets}--\ref{ss:description-cosets}.
 In fact, we have%
 \footnote{It is worth noting that applying the same construction to an element of 
$\ov S_\infty\times \ov S_\infty\times \ov S_\infty$ we generally
get a non-locally finite cell complex, which seems to
be pathological. In our case we get a disjoint union of finite complexes and all but
a finite number of components are double triangles with coinciding labels on white and black sides.}
 $$
 \K[\alpha]\setminus \G/\K[\beta]=  K[\alpha]\setminus G/K[\beta],
 $$
 therefore we get the same category $\cK$ of double cosets.
 For any unitary representation $\rho$ of $\G$ we have a $*$-representation of the category $\cK$
 (since a representation of $\G$ is also a representation of $G$).
 
 \sm
 
 {\bf \punct The equivalence theorem.}  See \cite{Ner-umn}, Sect.3.
 
 \begin{theorem}
  The map $\rho\mapsto \ov\rho$ establishes a bijection between the set of
  unitary representations of $\G$ and the set of $*$-representations of the category $\cK$.
 \end{theorem}
 
 \section{Spherical functions}
 
 \COUNTERS
 
 {\bf\punct Sphericity.}
 
 \begin{theorem}
  The pair $\G\supset \K$ is spherical,
  i.e., for any irreducible unitary representation of $\G$ the space of
  $\K$-fixed vectors has dimension $\le 1$.
 \end{theorem}
 
 {\sc Remark.}
 This is a trivial corollary of  multiplicativity Theorem \ref{th:multiplicativity} and the commutativity of
 the semigroup $\K\setminus\G/\K$. It is worth  noting that for the pair
 $$G_n=S_n\times S_n\times S_n,\qquad  K_n=\diag (S_n)$$
 the sphericity does not hold. \hfill $\boxtimes$
 
 \sm
 
  Consider a $\K$-fixed vector $v$ (a {\it spherical vector})
  for an irreducible representation $\rho$.
   Consider the following matrix element (a {\it spherical function}) 
 $$
 \Phi(p):=\la \rho(p)v,v\ra.
 $$
 It is easy to see that for any $r_1$, $r_2\in \K$ we have
 $\Phi(r_1 p r_2)=\Phi(p)$, i.e., $\Phi$ actually is a function
 on the double coset space $\K\setminus\G/\K$,
 i.e., on the space  $\wt \Xi[0,0]$ of checker surfaces.

\sm
 
 {\bf \punct Example.}
Here we present a construction of a family of representations of $\G$, for more constructions
 and a discussion of a classification problem, see \cite{Ner-umn}, Sect.3.
 Consider three Hilbert spaces $H^\blu$, $H^\red$, $H^\yel$ (they can be finite or infinite-dimensional).
 Consider the tensor product 
 $$
 X:=H^\blu \otimes H^\red\otimes H^\yel,
 $$
 and fix a unit vector $\xi\in X$.
 Consider an infinite tensor product%
 \footnote{Let $(Y_j,\upsilon_j)$ be a sequence of Hilbert spaces with distinguished unit vectors
 $\upsilon_j\in Y_j$. To define their tensor product $W$, we
 consider a system of formal vectors $y_1\otimes y_2\otimes\dots$ ({\it decomposable vectors})
 such that $y_j=\upsilon_j$ starting some  place. We assign for such vectors the following inner products
 \newline
 $\la y_1\otimes y_2\otimes\dots, \, y'_1\otimes y'_2\otimes\dots\ra=\prod_j \la y_j,y'_j\ra_{Y_j}$.
 \newline
 and assume that linear combinations of such vectors are dense in a Hilbert space $W$. This uniquely determines
 a Hilbert space $W$.  See \cite{Neu}, Section 4.1, or Addendum to \cite{Gui}.}
 \begin{equation}
 V:=(X,\xi)\otimes (X,\xi)\otimes(X,\xi)\otimes \dots
 \label{eq:otimes}
 \end{equation}
 Now let the first copy of $S_\infty$ acts by permutations of factors%
 \footnote{We emphasize that a finitely supported permutation of factors $H^\blue$ 
 send decomposable vectors to decomposable vectors,  for a permutation that is not finitely supported
 an image of a decomposable vector generally makes no sense.} $H^\blue$, the second copy by permutations
 of $H^\red$, and the third by permutations of $Y^\yel$. The subgroup $\K$ acts by permutations
 of factors $X$.

  \sm
 
 {\bf \punct An expression for spherical function.}
 Consider a representation $\rho$ of $\G$ constructed in the previous subsection.
 Consider the vector
 $$
 v:=\xi\otimes \xi\otimes \xi\otimes\dots \in V,
 $$
 it is a unique $\K$-invariant vector in
 $V$.
  We present a formula for the corresponding spherical function%
  \footnote{Generally speaking, $\rho(p)$ is reducible. However,
  the cyclic span $V'\subset V$ of $v$ is an irreducible representation of $\G$.
  Otherwise, we decompose $V'=W_1\oplus W_2$. Indeed, projections of $v$ to the both summands
  must be $\K$-invariant, but a $\K$-invariant vector is unique. Therefore $\Phi(p)$
  is the spherical function of the restriction of the representation $\rho(p)$
  to $V'$.
  \newline
  Assume that there are no proper subspaces $Z^\blu\subset H^\blu$, 
  $Z^\red\subset H^\red$, $Z^\yel\subset H^\yel$
  such that $\xi\in Z^\blu \otimes Z^\red\otimes Z^\yel$. Denote by
  $U$  the group consisting of all  triples of unitary operators
  ($A^\blu$, $A^\red$, $A^\yel$) acting in the corresponding spaces
  such that $A^\blu\otimes A^\red\otimes A^\yel \xi=\xi$. It is easy to show,
  that $U$ is a product of compact Lie groups and all but a finite number
  of factors are Abelian. For $v$ being in a general position 
  (and $\dim H^{\blu}$, $H^\red$, $H^\yel\ge 2$) the group $U$  is trivial. The group
  $U$ acts in each factor of the tensor product (\ref{eq:otimes}) 
  preserving  distinguished vector and therefore acts in the whole space $V$.
 There arises a natural conjecture: {\it $U$ and $\G$ are dual in the Schur--Weyl sense}
 (for  the spherical pair $S_\infty\times S_\infty\supset S_\infty$
 this property was obtained in \cite{Olsh-chip}, Subsect 5.5).
 If this is correct, then $\rho$ is irreducible for $\xi$ being in a general position.%
  }
 $$
 \Phi(p):=\la \rho(p)v,v\ra.
 $$
 Let us choose orthonormal bases
 $e_i^\blu\in H^\blu$, $e_j^\red\in H^\red$, $e_k^\yel\in H^\yel$,
 this defines a basis in the tensor product $X$.
We decompose $\xi\in X$ in this basis
 $$
 \xi=\sum_{i,j,k} \xi_{i,j,k}\, e_i^\blu\otimes e_j^\red\otimes  e_k^\yel\,.
 $$
 
 Fix a checker surface.
 Assign to each blue edge a vector $e_i^\blu$, to each red edge a vector
 $e^\red_j$, to each yellow edge a vector $e_k^\yellow$. Let us call such data  {\it assignments}.
 For a triangle $A$ denote by $i^\blu(A)$,  $j^\red(A)$, $k^\yel(A)$ the numbers of basis vectors
 on its sides.
 
 \begin{proposition} (see \cite{Ner-imrn})
  The spherical function $\Phi(\Xi)$ is given by
  \begin{multline*}
   \Phi(\Xi)=\sum\limits_{\text{\rm assignments}}\,\,
   \prod\limits_{\text{\rm white $A$}}
   \xi_{i^\blu(A),\,  j^\red(A),\, k^\yel(A)}
   \times \\ \times
     \prod\limits_{\text{\rm white $B$}}
  \ov{ \xi_{i^\blu(B),\,  j^\red(B),\, k^\yel(B)}}
  .
  \end{multline*}
 \end{proposition}

\section{Convolution algebras and their limit}

\COUNTERS

Here we explain how the multiplication of double cosets
arises as a limit of convolution algebras  on finite groups $G_n$.

\sm

{\bf\punct Convolution algebras.}
Let $H$ be a finite group,  $S$ its subgroup.
Denote by $\C[H]$ the group algebra of $H$ equipped with convolution
$*$.
Set 
$$
\delta_S:=\frac 1{\#(S)}\sum_{s\in S} \delta_h,
$$ 
where $\#(\cdot)$ denotes the number of elements in a finite set
and $\delta_h$ is the delta-function on $H$ supported by $h$.

Denote by $\C[S\setminus H/S]\subset \C[H]$ the space
of $S$-biinvariant functions, i.e., functions satisfying
the condition
$$
f(s_1 h s_2)=f(h), \text{where $s_1$, $s_2\in S$}.
$$
Equivalently,
\begin{equation}
f*\delta_S=\delta_S*f=f.
\label{eq:f-K}
\end{equation}
Clearly, $\C[S\setminus H/S]$ is a subalgebra in $\C[H]$.

Let $\rho$ be a unitary finite-dimensional representation of $H$ in the space $V$, 
denote by the same letter $\rho(\cdot)$ the corresponding representation of the group 
algebra. Denote by $V^S$ the space of $S$-fixed vectors in $V$.
Obviously, $\rho(\delta_S)$ is the operator of orthogonal projection to $V^S$.
For $f\in \C[S\setminus H/S]$ let us represent  $\rho(f)$ as 
a block  operator 
$$\rho(f): V^S\oplus (V^S)^\bot\to V^S\oplus (V^S)^\bot.$$
By (\ref{eq:f-K}), it has the form
$$
\rho(f)=\begin{pmatrix}
        \ov\rho(f)&0\\0&0
        \end{pmatrix}.
$$
Moreover operators $\ov\rho(f)$ determine a representation of the algebra 
$\C[S\setminus H/S]$ in $V^S$.
Standard arguments show that if $\rho$ is irreducible and $V^S\ne 0$,
then an initial representation $\rho$
is uniquely determined by representation $\ov\rho$.

There is a natural basis in $\C[S\setminus H/S]$. Namely, for each double coset
$\frp$, we consider the element of $\C[H]$ given by
$$
\delta_\frp:=\frac 1{\#\frp}\sum_{g\in\frp} \delta_g
.
$$
We have
\begin{equation}
\delta_\frp*\delta_\frq=\sum c_{\frp\frq}^\frr \delta_\frr,
\label{eq:structure-const}
\end{equation}
where $c_{\frp\frq}^\frr\in \Q $ are certain structure constants,
$$
\sum_\frr c_{\frp\frq}^\frr=1, \qquad  c_{\frp\frq}^\frr\ge 0
.
$$

\sm

{\sc Remark.} These considerations with minor variations can be repeated for a 
locally compact group $H$ and a compact subgroup $S$.
\hfill $\boxtimes$

\sm

{\sc Remark.}
An algebra of biinvariant functions (if we are able to understand it) is an important tool of
an investigation of representations of $H$. However, such algebras are difficult objects
and there few families of such objects that actually are explored.
The most important examples are: 

\sm

--- the Hecke--Iwahori algebras
(as $H=\GL$ over a finite field and $S$ is the group of all upper triangular matrices);

\sm

--- the
affine Hecke algebras (as $H=\GL$ over a $p$-adic field and $S$ is the Iwahori subgroup);

\sm

--- algebras of functions on a semisimple Lie groups
biinvariant with respect to its maximal compact subgroup (as $H=\GL(n,\R)$
and $S$ is the orthogonal group $O(n)$; or (a dual version) $H$ is a unitary group
$\U(n)$ and $S=\O(n)$).

\sm

Also recall that the Hecke and the affine Hecke algebras live by their own life 
independently on the initial groups. There is a huge literature on this topics.
\hfill $\boxtimes$

\sm

{\bf\punct Concentration of convolutions.}
Now consider a group 
$$H=G_n=S_n\times S_n\times S_n$$
 and 
the subgroup $K_n[\alpha]$ in the diagonal fixing 1, 2,\dots, $\alpha$.
Consider   algebras $\C\bigl[K_n[\alpha]\setminus G_n/K_n[\alpha]\bigr]$ for
fixed $\alpha$ and $n\to \infty$.

We assign a checker surface to an element of $K_n[\alpha]\setminus G_n/K_n[\alpha]$
as above.
 Thus we get an embedding
$$
\theta_n:
K_n[\alpha]\setminus G_n/K_n[\alpha]\,\longrightarrow \, \wt\Xi[\alpha,\alpha]
\simeq K[\alpha]\setminus G/K[\alpha].
$$
The image consists of all surfaces having $\le 2n$ triangles.
Let $\frp\in K[\alpha]\setminus G/K[\alpha]$ contain $2l$ triangles.
For each $n\in\N$ we define 
$\delta_{\frp(n)}\in \C\bigl[K_n[\alpha]\setminus G_n/K_n[\alpha]\bigr] 
$
by
$$
\delta_{\frp(n)}=\begin{cases}
                  0,\qquad&\text{if $n<l$};
                  \\
                  \theta_n^{-1} \frp ,\qquad&\text{if $n\ge l$}.
                 \end{cases}
$$

\begin{theorem}
Let $\frp$, $\frq$ be contained in the image of $\theta_m$,
i.e., they define elements of 
 $K_m[\alpha]\setminus G_m/ K_m[\alpha]$. Let $n>m$.
Decompose the following convolution as
$$
\delta_{\frp(n)}*\delta_{\frq(n)}=\sigma_n \delta_{(\frp\circledast\frq) (n)}+
\sum_{\frr\ne \frp\circledast\frq}b_{\frp(n)\frq(n)}^{\frr(n)} \delta_{\frr(n)}.
$$
Then the coefficient $\sigma_n$ tends to 1 as 
$n\to\infty$.
\end{theorem}

Thus the sum of remaining coefficients tends to 0. 

See a formal proof in \cite{Ner-concentration}, \cite{GN}, it is based on arguments
from \cite{Olsh-semigr}.

Multiplicativity Theorem \ref{th:multiplicativity} can be derived from 
the concentration of convolutions in a straightforward way, see \cite{GN}.

\section{A pre-limit convolution algebra}

\COUNTERS

{\bf \punct Algebras of conjugacy classes.}
Consider the group
$Q_n:=S_n\times S_n$, its diagonal $L_n\simeq S_n$ and the space of
conjugacy classes $Q_n/\!\!/S_n$.
As we noticed above,
$$
(S_n\times S_n)/\!\!/L_n\,\simeq\, K_n\setminus (S_n\times S_n\times S_n)/K_n\simeq \Xi_n
.$$
Consider the space $\C[Q_n/\!\!/L_n]$ of all $f\in \C[Q_n]$ satisfying
$$
f(h g h^{-1})=f(g), \qquad\text{where $g\in Q_n$, $h\in L_n$.}
$$
Obviously, this space is closed with respect to convolutions.
We intend to present a quasi-description of this algebra.

\sm

{\bf\punct The Ivanov-Kerov algebra.}
Consider a  disjoint union of all sets $Q_m/\!\!/L_m$
$$
\wh\Xi:=\coprod_{m=0}^\infty \Xi_n\,\simeq \,\coprod_{m=0}^\infty Q_m/\!\!/L_m,
$$
({\it below use the same notation for a conjugacy class
and the corresponding checker surface}).

\sm

{\sc Remark.} Emphasize that $\wh \Xi\ne \Xi[0,0]$. Indeed,
an element of $\wh \Xi$ is an arbitrary finite
checker surface, and an element of $\Xi[0,0]$ is a checker surface without double triangles.
A removing double triangles determines a map $\wh\Xi\to \Xi[0,0]$.
A preimage of any point is countable.
\hfill $\boxtimes$

\sm

For any $n$ and $\frp\in Q_n/\!\!/L_n$ we define
an element $\Delta_{\frp}^{n,n} $ of the group algebra $\C[Q_n]$ by
$$
\Delta_{\frp}^{n,n} :=\sum_{g\in \frp} \delta_g.
$$
Let $m\ge n$ and
$\frp\in Q_n/\!\!/L_n$.
Consider the corresponding element  $\frp(m)\in Q_m/\!\!/L_m$,
i.e., we take the corresponding checker surface and add $m-n$
double triangles. We set 
$$
\Delta_{\frp}^{m,n}:=\frac{(N-n)!}{(N-m)!} \Delta_{\frp}^{n,n}.
$$

\begin{theorem}
There exists  an associative  algebra $\frA$ with a basis $u_\frp$, where $\frp$ ranges in $\wh\Xi$,
and 
a product
$$
u_{\frp}\circ u_\frq=\sum_{\frq\in \wh\Xi}
a_{\frp\frq}^\frr u_\frr
$$
such that

\sm

$\bullet$ $a_{\frp\frq}^\frr\in\Z_+$.

\sm

$\bullet$ Consider the linear map $\Pi_n: \frA\to \C[Q_n/\!\!/L_n]$
defined on $u_\frp$, $\fru\in Q_k/\!\!/L_k$ by
$$
\Pi_n u_\frp:=
\begin{cases}
 \Delta_{\frp}^{n,k}, \qquad& \text{if $k\le n$}
 \\
 0 , \qquad &\text{if $k>n$}.
\end{cases}
$$
Then $\Pi_n$ is 
 a homomorphism from $\frA$ to $\C[Q_n/\!\!/L_n]$.

\sm

$\bullet$ Let $\frp\in Q_n/\!\!/L_n$, $\frq\in Q_m/\!\!/L_m$, $\frr\in Q_l/\!\!/L_l$. 
Then $a_{\frp\frq}^\frr$ are zero if $l>m+n$ or if $l<\min(m,n)$. 
\end{theorem}

A similar statement for $\C[S_n/\!\!/S_n]$ was obtained by Ivanov and Kerov
in \cite{IK}, it was extended to more general objects including
$\C[Q_n/\!\!/L_n]$ in \cite{Ner-IK}.

\sm

{\bf\punct A description of the algebra $\frA$.}
See \cite{Ner-IK}.
Let $X$, $Y$ be finite sets. We say that a {\it partial bijection}
$s:X\to Y$ is a bijection of a subset $B\subset X$ to a subset $B\subset Y$,
we say $\dom(s):=A$, $\im (s):=B$.

Let $\frp$, $\frq\in \wh \Xi$. Consider a partial bijection
$s$ from the set of black triangles of $\frp$ to the set of white triangles 
of $\frq$. For each $A\in \dom s$ we remove the interior of the black triangle $A\subset \frp$ 
and the interior of the white triangle
$s(A)\in\frq$ and identify their boundaries  according colors of  edges.
We normalize the resulting cell complex and come to
 a new checker surface $\frp\odot_{s} \frq$.
The formula for the product in $\frA$ is
$$
u_\frp\circ u_\frq=\sum_s u_{\frp \odot_{s} \frq},
$$
where the summation is taken over all partial bijection.

\sm

{\bf \punct The Poisson algebra.} Denote by $\frA_n$ the subspace in
$\frA$ generated by all $u_\frp$, where $\frp$ ranges in $\coprod_{i=0}^n Q_n/\!\!/K_n$.
The get  a filtration,
$$
v\in\frA_n,\,\, w\in \frA_m \quad\Rightarrow \quad v\circ w\in \frA_{m+n}.
$$
The corresponding graded associative algebra 
$\oplus_{j=0}^\infty \frA_{n}/\frA_{n-1}$
is
commutative, a product of basis elements is given by
$$
u_\frp\cdot u_\frq=u_{\frp\coprod \frq},
$$
where $\frp\coprod\frq$ is a disjoint union of checker surfaces.

Also, we  have a structure of a  graded Poisson Lie algebra on $\oplus_{n=0}^\infty \frA_{n}/\frA_{n-1}$,
namely, for any $x\in \frA_m$, $y\in\frA_n$ we define their bracket $\{x,y\}$
as the image of $x\ast y-y\ast x$ in $\frA_{m+n-1}/\frA_{m+n-2}$.
Let us describe the bracket. Consider set $M(\frp,\frq)$ of all pairs $\phi=(B,C)$, where
$B$ is a black triangle of $\frp$ and $C$ is a white triangle of $\frq$. For
any $\psi$ we define a surface $\frp\circledcirc_\phi\frq$, it is obtained
by removing the interior of $B$ from $\frp$ and the interior of $C$ from $\frq$ and identification 
of boundaries of these triangles. 
In this notation, 
$$
\{u_\frp, u_\frq\}=
\sum_{\phi\in M(\frp,\frq)} u_{\frp\circledcirc_\phi \frq} -
\sum_{\psi\in M(\frq,\frp)}u_{\frq\circledcirc_\psi \frp}.
$$

	\noindent
	\tt Math.Dept., University of Vienna,
	\\
	Oskar-Morgenstern-Platz 1, 1090 Wien;
	\\
	\& Institute for Theoretical and Experimental Physics (Moscow);
	\\
	\& Mech.Math.Dept., Moscow State University;
	\\
\&	Institute for information transmission problems (Moscow);
	\\
	e-mail: yurii.neretin(at) univie.ac.at
	\\
	URL:www.mat.univie.ac.at/$\sim$neretin

\end{document}